\documentclass[paper=letter, fontsize=12pt]{article} 
\RequirePackage[T1]{fontenc}
\RequirePackage{amsthm,amsmath,amsfonts}
\usepackage{longtable}
\usepackage{graphicx}
\usepackage{mathtools}
\usepackage{latexsym}
\usepackage{multirow}
\usepackage{setspace}
\usepackage{adjustbox}
\usepackage[bf]{caption}
\usepackage{emptypage}
\usepackage{subfigure}
\usepackage{float}
\usepackage{algorithm}
\usepackage{algpseudocode}
\usepackage{booktabs}
\usepackage{siunitx}

\RequirePackage[round]{natbib}
\RequirePackage[colorlinks,citecolor=blue,urlcolor=blue]{hyperref}

\allowdisplaybreaks[1]

\usepackage{times}
\usepackage{keyval}
\usepackage{calc}

\usepackage[letterpaper]{geometry} 
\renewcommand{\baselinestretch}{1.1}

\usepackage{fancyhdr}
\newcommand{\iu}{\mathrm{i}\mkern1mu}

\numberwithin{equation}{section}
\theoremstyle{plain}
\newtheorem{theorem}{Theorem}[section]
\newtheorem{assumption}{Assumption}[section]
\newtheorem{remark}{Remark}[section]
\newtheorem{lemma}{Lemma}[section]

\newtheorem{critere}{Criterion}[theorem]

\newtheorem{proposition}{Proposition}[section]
\newtheorem{variational inequality}{Variational inequality}[section]

\DeclareSymbolFont{symbols}{OMS}{cmsy}{m}{n}
\DeclareFontFamily{OT1}{pzc}{}
\DeclareFontShape{OT1}{pzc}{m}{it}{<-> s * [1.10] pzcmi7t}{}
\DeclareMathAlphabet{\mathbf}{OT1}{cmr}{bx}{n}
\DeclareMathAlphabet{\mathsf}{OT1}{cmss}{m}{n}
\DeclareMathAlphabet{\mathit}{OT1}{cm}{m}{it}
\DeclareMathAlphabet{\mathpzc}{OT1}{pzc}{m}{it}
\DeclareMathAlphabet{\pazocal}{OMS}{zplm}{m}{n}

\DeclareGraphicsExtensions{.eps}
\allowdisplaybreaks[1]

\usepackage{fancyhdr}
\renewcommand{\dot}[1] {\overset{\,_{\mbox{\Large .}}}{#1}}
\renewcommand{\ddot}[1] {\overset{\,_{\mbox{\Large ..}}}{#1}}

\DeclareSymbolFont{boldoperators}{OT1}{cmr}{bx}{n}
\SetSymbolFont{boldoperators}{bold}{OT1}{cmr}{bx}{n}
\edef\bar{\unexpanded{\protect\mathaccentV{bar}}\number\symboldoperators16}

\begin{document}
	\rhead{\textit{\today}}
	\lhead{\textit{Gao \& Hyndman}}
	\chead{Boundary error control for BSDE-CFFT}

	\title{Boundary error control for numerical solution of BSDEs by the convolution-FFT method}
	
	\author{
		Xiang Gao\footnote{ 
			Department of Mathematics and Statistics, 
			Concordia University, 
			1455 Boulevard de Maisonneuve Ouest,
			Montr\'eal, Qu\'ebec,
			Canada H3G 1M8.
		}
		\ and 
		Cody Hyndman\footnotemark[1] \, \footnote{Corresponding Author: cody.hyndman@concordia.ca}
	}
	
	\date{November 3, 2025}
	
	\maketitle
	
	\abstract{We first review the convolution fast-Fourier-transform (CFFT) approach for the numerical solution of backward stochastic differential equations (BSDEs) introduced in \citep{hyndman2017convolution}. We then propose a method for improving the boundary errors obtained when valuing options using this approach. We modify the damping and shifting schemes used in the original formulation, which transforms the target function into a bounded periodic function so that Fourier transforms can be applied successfully.  Time-dependent shifting reduces boundary error significantly. We present numerical results for our implementation and provide a detailed error analysis showing the improved accuracy and convergence of the modified convolution method.}
        
	\vspace{5mm}
	

	\vspace{5mm}

	\renewcommand{\baselinestretch}{1.5}
		
	\section{Introduction\label{sec:Intro}}
In this paper we study a convolution-based numerical method for solving backward stochastic differential equations (BSDEs) that arise in option valuation. Specifically, we consider the coupled forward--backward system
\begin{align}
X_t &= x + \int_0^t \eta(s,X_s)\,ds + \int_0^t \sigma(s,X_s)\,dW_s,  \label{FSDE}\\
Y_t &= g(X_T) + \int_t^T f(s,X_s,Y_s,Z_s)\,ds - \int_t^T Z_s^{\top}\,dW_s, \label{BSDE}
\end{align}
on a filtered probability space $(\Omega,\mathcal{F},\{\mathcal{F}_t\}_{t\in[0,T]},\mathbb{P})$, where $W$ is an $n$-dimensional Brownian motion and $\{ \mathcal{F}_t \}$ is its augmented natural filtration. Given the coefficient functions $\eta$, $\sigma$, $f$, and terminal condition $g$, the unknown adapted processes are $(X_t,Y_t,Z_t)$.

 \cite{pardoux1990adapted} first introduced nonlinear BSDEs and established foundational existence and uniqueness results. Subsequent work developed the theory in several directions, including forward--backward SDEs (FBSDEs) \cite{antonelli1993backward}, solvability and well-posedness results \cite{yong1997finding}, and extensions beyond the globally Lipschitz setting, including quadratic-growth drivers \cite{kobylanski2000backward}.

BSDEs have many applications in mathematical finance, including pricing, hedging, and nonlinear valuation problems. As a result, numerical methods for BSDEs have been studied extensively. Common approaches include simulation-based schemes \cite{bouchard2004discrete}, PDE-based discretizations \citet{douglas1996numerical}, Picard iteration \cite{bender2007forward}, and regression-based methods \cite{lemor2006rate}. Spatial discretization methods for Markovian BSDEs can be found, for example, in \cite{huijskens2016efficient} and in \cite{hyndman2017convolution}.

	The convolution method for numerical solution of BSDEs introduced in  \citet{hyndman2017convolution} is fast and accurate, and was further studied in \cite{polynice2022convolution} through alternative spatial approximation schemes aimed at reducing extrapolation error. However, in the original formulation, numerical errors can grow rapidly as the initial value of the forward process approaches the truncation boundaries. This paper is complementary to our convolution--FFT methods for option pricing based on characteristic functions under the Heston model \cite{wang_hynd:heston}. While the Heston model convolution--FFT approach focuses on Fourier inversion of terminal payoffs, the present paper develops a convolution-based time-stepping scheme for general BSDEs.  Throughout, the convolution representation is understood as a short-time approximation consistent with the underlying time discretization of the BSDE.

The paper is organized as follows. In Section~\ref{sec:review}, we review the convolution method with an implicit Euler time discretization scheme for BSDEs. In Section~\ref{sec:improve}, we introduce the proposed boundary-error control method and provide an error analysis describing how local errors vary with respect to the initial value. In Section~\ref{sec:application}, we present numerical results for option valuation and compare the performance of the modified method with the original scheme of \cite{hyndman2017convolution}. Section~\ref{sec:conclusion} concludes and an appendix contains proofs of technical results.

	\section{Assumptions and the convolution method}\label{sec:review}
	In this section, we will review the convolution approach to the numerical solution of  BSDEs introduced in \citet{hyndman2017convolution}.  To ensure the existence and uniqueness of an adapted solution to (\ref{FSDE})-(\ref{BSDE}), we require the following conditions (see e.g. \citet{ma1999forward}) to be satisfied.
        \begin{assumption}\label{assumption1}
          For  $\eta :[0,T]\times\mathbb{R} \rightarrow \mathbb{R}$, $\sigma : [0,T]\times\mathbb{R}\rightarrow\mathbb{R}$, $f :[0,T]\times\mathbb{R}\times\mathbb{R}\times\mathbb{R} \rightarrow\mathbb{R}$ and $g: \mathbb{R} \rightarrow \mathbb{R}$ assume:
\begin{itemize}
\item[(i)] The functions $\eta$, $\sigma$, $f$, and $g$ are uniformly Lipschitz continuous with bounded first order derivatives in the space variables, for all $t\in[0,T]$:
	\begin{align*}
	|\eta(t,u) - \eta(t,v)| \leq&  K(|u - v|),\\
	|\sigma(t,u) - \sigma(t,v)| \leq&K(|u - v|),\\
	|g(u) - g(v)|  \leq&  K(|u - v|),\\
	|f(t,\xi) - f(t,\zeta)| \leq& K(\|\xi - \zeta\|_{\infty}),
	\end{align*}
for some constant $K$ independent of $u,v\in\mathbb{R}$ and $\xi,\zeta\in\mathbb{R}\times\mathbb{R}\times\mathbb{R}$.
			\item[(ii)] The volatility $\Sigma(t,x) = {\sigma(t,x)} {\sigma}^\top(t,x)$ is continuous and $L^{2}$-bounded 
			$
				\|\Sigma\|_2\leq C,
			$
			for some positive constant $C$.
		\end{itemize}
	\end{assumption}

        Using a uniform time discretization with step size $\Delta t = t_{i+1}-t_i$ we write 
        $X_{i} = X_{t_{i}}$, $Y_{i} = Y_{t_{i}}$, $Z_{i}=Z_{t_{i}}$, and $\Delta W_{i} = (W_{t_{i+1}}-W_{i})$.  Applying a first-order Euler  scheme we have
$$X_{i+1} = X_{i} + \eta(t_{i},X_{i}) \Delta t + \sigma(t_{i},X_{i})\Delta W_{i}; \qquad X_{0}=x; \qquad i=0,\ldots,N-1$$
for the approximation of the forward SDE (\ref{FSDE}) and 
 $$Y_{{i}} = Y_{{i+1}} +  f(i,X_{i},Y_{i},Z_{i})\Delta t - Z_{i} \Delta W_{i}; \qquad Y_{N}=g(X_{N}), \qquad Z_{N}=0 \qquad i=N-1,\ldots, 0$$
for the implicit approximation of the backward SDE (\ref{BSDE}).  
By taking the conditional expectation of $Y_{i}$ and $Y_{i}\Delta W_{i}$ given $X_{i}$ we obtain
the standard explicit scheme
\begin{align}
Y_i &= \dot{Y}_{i} + \Delta t\, f\!\left(t_i,X_{i}, \mathbb{E}[Y_{i+1}\mid X_{i}], Z_i\right),  \label{euler1}   \\
Z_i &= \frac{1}{\Delta t}\,\mathbb{E}\!\left[\, {Y}_{i+1}\Delta W_i\,\middle|\,X_{t_i}\,\right], \label{euler2} \\
\dot{Y}_{i} &= \mathbb{E}\!\left[Y_{i+1} \middle| \,X_{i}\,\right]. \label{euler3}
\end{align}
Provided the conditional expectations can be calculated, at least approximately, we have a recursive numerical method for solving (\ref{FSDE})-(\ref{BSDE}).  The $Z$-component is approximated consistently with the spatial discretization underlying the convolution scheme.  %

Under standard regularity assumptions on the coefficients of the forward SDE and the driver of the BSDE, the solution admits a Markovian representation. In particular, there exists a deterministic function $u:[0,T]\times\mathbb{R}\to\mathbb{R}$ such that
\begin{align}
  Y_t &= u(t,X_t),\label{markovian_repY} \\
  Z_{t} &= \sigma(t,X_{t})^\top \nabla_{x}u(t,X_{t}) \label{markovian_repZ}
\end{align}
$t\in [0,T]$, almost surely. The representation (\ref{markovian_repY})-(\ref{markovian_repZ}) follows from the nonlinear Feynman--Kac formula, see \cite{pardoux1990adapted}, and allows the backward component of the BSDE to be expressed as a function of time and the current value of the forward process.

Conditioning on $X_{t_i}=x$ and denoting $\Delta t = t_{i+1}-t_i$, the conditional expectations appearing in the time--discretized BSDE may be written in integral form as
\begin{equation}\label{cond_exp_integral}
\dot{Y}_{i}(x) = \mathbb{E}\!\left[ Y_{t_{i+1}} \mid X_{t_i}=x \right]
=
\int_{\mathbb{R}} u(t_{i+1},x') \, \phi(x' \mid t_i,x)\, dx',
\end{equation}
where $\phi(\cdot \mid t_i,x)$ denotes the transition density of the forward process $X_{t_{i+1}}$ conditional on $X_{t_i}=x$.

For small time increments $\Delta t$, the transition density admits a short--time approximation consistent with the local behavior of the forward diffusion. In particular, conditional on $X_{t_i}=x$, the density $\phi(x' \mid t_i,x)$ may be approximated by a Gaussian density of the form
\begin{equation}\label{short_time_density}
\phi(x' \mid t_i,x)
\approx
\phi(x'-x \mid t_i,x)
=
\frac{1}{\sqrt{2\pi \sigma^2(t_n,x)\Delta t}}
\exp\!\left(
-\frac{\big(x' - x - \eta(t_i,x)\Delta t\big)^2}{2\sigma^2(t_i,x)\Delta t}
\right),
\end{equation}
where $\eta(t,x)$ and $\sigma(t,x)$ are the drift and diffusion coefficients of the forward SDE.
Such short–time Gaussian approximations of the transition density are standard for diffusion processes and are consistent with Euler–Maruyama discretization; see, for example, \citet{risken1996fokker} or classical results on transition densities of SDEs.

Substituting \eqref{short_time_density} into \eqref{cond_exp_integral}, the conditional expectation can be written in convolution form,
\begin{equation}\label{convolution_form}
\dot{Y}_{i}(x) = \mathbb{E}\!\left[ Y_{t_{i+1}} \mid X_{t_i}=x \right]
\approx
\int_{\mathbb{R}} u(t_{i+1},x') \, \phi(x'-x \mid t_i,x)\, dx'
=
\big( u(t_{i+1},\cdot) * \phi(\cdot \mid t_i,x) \big)(x),
\end{equation}
where $*$ denotes convolution with respect to the spatial variable.

This short–time convolution representation provides the basis for the numerical scheme developed below. In practice, the convolution integral is evaluated on a truncated spatial domain and discretized using FFT-based techniques, leading to an efficient approximation of the conditional expectations appearing in the BSDE time–stepping scheme.

A similar integral representation is available for the conditional expectation appearing in the standard discrete approximation of the $Z$--component. In particular, for the Euler-type scheme
\begin{equation}\label{eulerZ}
Z_i(x) = \frac{1}{\Delta t}\,\mathbb{E}\!\left[\, Y_{i+1}\Delta W_i\,\middle|\,X_{t_i}=x\,\right],
\end{equation}
we may write, conditioning on $X_{t_i}=x$ and using the Markovian representation $Y_{i+1}=u(t_{i+1},X_{t_{i+1}})$,
\begin{equation}\label{Z_integral_1}
\mathbb{E}\!\left[\, Y_{i+1}\Delta W_i\,\middle|\,X_{t_i}=x\,\right]
=
\int_{\mathbb{R}} u(t_{i+1},x') \,
\mathbb{E}\!\left[\,\Delta W_i \,\middle|\, X_{t_{i+1}}=x',\,X_{t_i}=x\,\right]\,
\phi(x' \mid t_i,x)\,dx',
\end{equation}
where $\phi(\cdot \mid t_i,x)$ denotes the transition density of $X_{t_{i+1}}$ conditional on $X_{t_i}=x$.

Under the short--time Euler approximation of the forward diffusion,
\begin{equation}\label{euler_forward_local}
X_{t_{i+1}} \approx x + \eta(t_i,x)\Delta t + \sigma(t_i,x)\Delta W_i,
\end{equation}
we have the identity
\begin{equation}\label{dW_as_function_of_increment}
\Delta W_i \approx \frac{X_{t_{i+1}}-x-\eta(t_i,x)\Delta t}{\sigma(t_i,x)}.
\end{equation}
Substituting \eqref{dW_as_function_of_increment} into \eqref{Z_integral_1} and using the corresponding short--time approximation $\phi(x'\mid t_i,x)\approx \phi(x'-x\mid t_i,x)$ yields
\begin{equation}\label{Z_integral_2}
\mathbb{E}\!\left[\, Y_{i+1}\Delta W_i\,\middle|\,X_{t_i}=x\,\right]
\approx
\frac{1}{\sigma(t_i,x)}\int_{\mathbb{R}} u(t_{i+1},x')\,
\big(x'-x-\eta(t_i,x)\Delta t\big)\,
\phi(x'-x \mid t_i,x)\,dx'.
\end{equation}
Therefore,
\begin{equation}\label{Z_convolution_form}
Z_i(x)
\approx
\frac{1}{\Delta t\,\sigma(t_i,x)}
\int_{\mathbb{R}} u(t_{i+1},x')\,
\big(x'-x-\eta(t_i,x)\Delta t\big)\,
\phi(x'-x \mid t_i,x)\,dx'.
\end{equation}
Defining the kernel
\[
\kappa_i(z;x) := \frac{z-\eta(t_i,x)\Delta t}{\Delta t\,\sigma(t_i,x)}\,\phi(z \mid t_i,x),
\qquad z\in\mathbb{R},
\]
the approximation \eqref{Z_convolution_form} can be written compactly as the convolution
\begin{equation}\label{Z_convolution_compact}
Z_i(x) \approx \big(u(t_{i+1},\cdot) * \kappa_i(z;\cdot)\big)(x),
\end{equation}
where $*$ denotes convolution with respect to the spatial variable.

The convolution representations derived above provide a convenient framework for numerical approximation. In particular, convolution integrals can be evaluated efficiently in the Fourier domain using the convolution theorem, which transforms convolutions in physical space into pointwise products in frequency space. At the continuous level, this amounts to applying the Fourier transform to the convolution kernels and the functions being convolved, followed by inversion of the transform to recover the approximations in physical space.  Further details can be found in \citet{hyndman2017convolution} and \cite{polynice2022convolution}.

In practice, the spatial domain must be truncated and the resulting integrals discretized. The numerical evaluation of these Fourier-domain expressions is carried out using discrete Fourier transforms and fast Fourier transform (FFT) algorithms. Since truncation, discretization, and boundary effects play a crucial role in the stability and accuracy of the method, we defer the detailed discussion of Fourier discretization, damping and shifting techniques, and FFT-based implementation to the next section.

\section{Boundary control schemes}\label{sec:improve}
This section develops boundary-control modifications for the convolution-based FBSDE scheme. We first motivate damping and shifting transformations that reduce truncation and periodic-extension artifacts, then derive the corresponding Fourier-domain representations and present an FFT-based implementation together with error estimates.

\subsection{Boundary effects and damping--shifting strategy}\label{sec3.1}

In practice, the Fourier transform is applied on a truncated spatial domain. While truncation allows the transform to be computed numerically, it also introduces boundary effects due to the implicit periodic extension associated with Fourier-based methods. If not properly controlled, these boundary effects can lead to large numerical errors and instability in the resulting approximations.

Convergence results for Euler discretizations of BSDEs are available under suitable regularity conditions. For example, \citet{zhang2004numerical} show that, for a partition scheme $\Delta$, the approximation error satisfies
\[
\mathbb{E}\left|Y(X)-Y(X^\Delta)\right|^2 \leq C\left(1 + \left|x\right|^2\right)\left|\Delta\right|,
\]
provided the target function satisfies an appropriate Lipschitz condition. However, in option pricing problems, terminal payoff functions are typically non-Lipschitz and unbounded, and additional modifications are required to ensure numerical stability.

To address the lack of integrability, \citet{carr1999option} introduced an exponential damping factor to enforce decay of the target function prior to applying the Fourier transform. Related damping ideas have been adopted in \citet{lord2007optimal} and \citet{hyndman2017convolution}, where negative damping parameters are used to improve stability. While damping ensures integrability, truncation alone does not eliminate boundary artifacts. In particular, the implicit periodic extension induced by the Fourier transform may still generate significant errors near the boundaries of the computational domain.

To mitigate these effects, \citet{hyndman2017convolution} proposed a linear shifting function to enforce periodicity of the target function prior to applying the Fourier transform. However, linear shifting may itself introduce large boundary errors when applied to option payoff functions. This motivates the use of a shifting function adapted to the structure of the terminal condition.

In this work, we combine exponential damping with an exponential shifting function tailored to option pricing payoffs. Specifically, we introduce a modified target function
\begin{equation}\label{shift2}
\tilde u(x) = e^{\alpha x}\big(u(x) - h(x)\big),
\end{equation}
where $\alpha<0$ is a damping parameter and $h(x)=Ae^{x}+B$ is chosen so that the modified function $\tilde u$ satisfies the periodicity conditions
\[
\tilde u(x_0) = \tilde u(x_N), \qquad
\frac{d\tilde u}{dx}(x_0) = \frac{d\tilde u}{dx}(x_N),
\]
on the truncated spatial interval $[x_0,x_N]$. This construction significantly reduces boundary-induced errors and improves numerical stability.

With this damped and shifted formulation, the convolution representations derived in Section~2 can be evaluated efficiently in the Fourier domain using the convolution theorem. At the continuous level, this leads to Fourier-domain expressions for the approximations of the $Y$- and $Z$-components, which can subsequently be inverted to recover the solutions in physical space. The precise Fourier transform definitions, damping and shifting in frequency space, and discrete FFT-based implementation are presented in the following subsections.

\subsection{Fourier-domain representation of the damped--shifted scheme}\label{sec3.2}

With the damped and shifted formulation introduced in Section~\ref{sec3.1}, the convolution updates for the backward components can be expressed conveniently in the Fourier domain. At the continuous level, convolution in physical space corresponds to pointwise multiplication in frequency space, and the original target functions can be recovered by applying the inverse Fourier transform.

Under the short-time Gaussian approximation over a single step $\Delta t$, the forward increment satisfies
\[
X_{t+\Delta t}-X_t \approx \eta\,\Delta t + \sigma\,\Delta W_t,
\]
so its characteristic function is
\begin{equation}\label{psi_def}
\psi(v)
:= \mathbb{E}\!\left[e^{\iu v (X_{t+\Delta t}-X_t)}\mid X_t\right]
= \exp\!\left(\Delta t\left(\eta\,\iu v - \tfrac{1}{2}\sigma^2 v^2\right)\right),
\end{equation}
with derivative $\psi'(v)=\Delta t(\eta\,\iu-\sigma^2 v)\psi(v)$.
In the damped formulation \eqref{shift2}, the Fourier multipliers corresponding to the $Y$- and $Z$-updates are evaluated at complex-shifted frequencies. 

We first consider the update for the $Y$--component. Applyin the damping and shifting transformation \eqref{shift2} to \eqref{convolution_form}, taking the Fourier transform, applying the convolution theorem, inverting the Fourier transform, then inverting the damping and shifting gives that the approximation $\dot Y_t$ can be written as
\begin{align}
\dot Y_t(x)
= {} & \mathfrak{F}^{-1}\!\left[
\mathfrak{F}\!\left[ Y_{t+1} \right](v)\,\Psi_y(v)
\right](x) \notag\\
= {} & e^{-\alpha x}\,
\mathfrak{F}^{-1}\!\left[
\mathfrak{F}\!\left[ \tilde Y_{t+1}^{(\alpha)} \right](v)\,\Psi_y(v)
\right](x)
+ A\,\mathbb{E}\!\left[ e^{X_{t+1}} \mid X_t=x \right] + B
\label{Y_recovery}
\end{align}
where we define
\begin{equation}\label{Psi_def}
\Psi_y(v) := \psi(v+\alpha\iu),
\end{equation}
for the fixed damping parameter $\alpha<0$.
Under the short--time Gaussian approximation of the forward process \eqref{euler_forward_local}, the conditional expectation in \eqref{Y_recovery} can be evaluated explicitly, yielding
\begin{equation}
\dot Y_t(x)
=
e^{-\alpha x}\,
\mathfrak{F}^{-1}\!\left[
\mathfrak{F}\!\left[ \tilde Y_{t+1}^{(\alpha)} \right](v)\,\Psi_y(v)
\right](x)
+ A e^{x}\psi(-\iu) + B.
\end{equation}
 The explicit recovery terms associated with the exponential shift $h(x)=Ae^x+B$ are given by
\begin{equation}\label{ddot_def}
\ddot Y(x) := e^x\psi(-\iu),\qquad
\ddot Z(x) := -\,\frac{\eta\,\Delta t\,e^x\psi(-\iu)+\iu\,e^x\psi'(-\iu)}{\sigma\,\Delta t}.
\end{equation}
These recovery terms can be precomputed on the spatial grid.

A similar representation holds for the $Z$--component. Applying the damping and shifting transformation \eqref{shift2} to the approximation \eqref{Z_convolution_compact} for $Z_t$, taking the Fourier transform, applying the convolution theorem, inverting the Fourier transform and inverting the damping and shifting, we obtain
\begin{align}
Z_t(x)
= {} & \mathfrak{F}^{-1}\!\left[
\mathfrak{F}\!\left[ Y_{t+1} \right](v)\,\Psi_z(v)
\right](x) \notag\\
= {} & e^{-\alpha x}\,
\mathfrak{F}^{-1}\!\left[
\mathfrak{F}\!\left[ \tilde Y_{t+1}^{(\alpha)} \right](v)\,\Psi_z(v)
\right](x)
+ \frac{A}{\Delta t}\,
\mathbb{E}\!\left[ e^{X_{t+1}}\Delta W_t \mid X_t=x \right]
\label{Z_recovery}
\end{align}
where we define
$$\Psi_z(v) := \sigma(\iu v-\alpha)\,\psi(v+\alpha\iu).$$
Evaluating the conditional expectation using the short--time approximation of the forward increment \eqref{euler_forward_local} gives
\begin{equation}
Z_t(x)
=
e^{-\alpha x}\,
\mathfrak{F}^{-1}\!\left[
\mathfrak{F}\!\left[ \tilde Y_{t+1}^{(\alpha)} \right](v)\,\Psi_z(v)
\right](x)
-
A\,\left(\frac{\eta\,\Delta t\, e^{x}\psi(-\iu) + \iu\, e^{x}\psi'(-\iu)}{\sigma\,\Delta t}\right).
\end{equation}
Here,
\[
\psi(v) = \exp\!\left(
\Delta t\left( \eta\,\iu v - \tfrac{1}{2}\sigma^2 v^2 \right)
\right)
\]
denotes the characteristic function associated with the short--time Gaussian approximation of the forward increment, and
\[
\psi'(v) = \Delta t\left( \eta\,\iu - \sigma^2 v \right)\psi(v)
\]
is its first derivative with respect to $v$.

\subsection{Discrete Fourier implementation}\label{sec3.3}

We now describe the discrete Fourier implementation of the damped--shifted convolution scheme. All Fourier-domain formulas in Section~\ref{sec3.2} are evaluated numerically on a truncated spatial grid and the corresponding frequency grid. Let the truncation interval be $[x_0,x_N]$ with length $L=x_N-x_0$, and define the uniform spatial grid
\[
x_n = x_0 + n\Delta x,\qquad n=0,1,\dots,N-1,\qquad \Delta x=\frac{L}{N}.
\]
The associated frequency grid is
\[
v_k = \left(k-\frac{N}{2}\right)\Delta v,\qquad k=0,1,\dots,N-1,\qquad \Delta v=\frac{2\pi}{L},
\]
so that the Nyquist relation $\Delta x\,\Delta v=2\pi/N$ holds.

For a grid function $\{f_n\}_{n=0}^{N-1}$, we write $\mathfrak{D}[f](v_k)$ for its discrete Fourier transform and $\mathfrak{D}^{-1}$ for the inverse transform, defined by
\begin{align}
\mathfrak{D}[f](v_k)
&:= \sum_{n=0}^{N-1} f_n\,e^{-\iu v_k x_n}, \label{DFT_def}\\
\mathfrak{D}^{-1}[F](x_n)
&:= \frac{1}{N}\sum_{k=0}^{N-1} F(v_k)\,e^{\iu v_k x_n}. \label{IDFT_def}
\end{align}
In order to work with a frequency grid centered at zero, we use the standard phase-shift identity
\[
e^{\pm \iu v_k x_n}
= (-1)^{n+k}\,e^{\pm \iu (k\Delta v)(n\Delta x)},
\]
which results in the $(-1)^n$ factors appearing in the FFT-ready formulas below. In practice, this centering is implemented by multiplying the input and output grid vectors componentwise by $(-1)^n$.

We summarize the resulting FFT-based backward iteration with boundary control in Algorithm~\ref{Algo2}. The shifting parameters $(A,B)$ are updated at each time step to enforce periodicity of the modified function $\tilde u$, while the damping parameter $\alpha$ is fixed throughout. This choice reduces computational overhead and avoids the stability issues that can arise when $\alpha$ is updated adaptively.

        \begin{algorithm}[t]
\caption{Convolution--FFT scheme for the FBSDE system with damping and shifting}
\label{Algo2}
\begin{algorithmic}[2]
\Require Truncation length $L$, grid size $N$, time grid $\{t_k\}_{k=0}^n$ with $\Delta t = T/n$, damping parameter $\alpha<0$ (fixed),
        payoff/terminal function $g$, driver $f$, model parameters $(\eta,\sigma)$.
\Ensure Approximations $\{\bar Y_k^\Delta(x_i), \bar Z_k^\Delta(x_i)\}$ on the spatial grid.

\State Set $\Delta x \gets L/N$, $\Delta v \gets 2\pi/L$
\State Define spatial grid $x_i \gets x_0 + i\Delta x$, $i=0,\dots,N-1$
\State Define frequency grid $v_j \gets (j-\frac{N}{2})\Delta v$, $j=0,\dots,N-1$

\State Initialize terminal values: $Y_{n,i}^\Delta \gets g(x_i)$ for $i=0,\dots,N-1$
\State Initialize $Z_{n,i}^\Delta$ (e.g.\ regression/finite-difference/analytic if available) for $i=0,\dots,N-1$

\State Precompute Fourier multipliers for all $j=0,\dots,N-1$:
\[
\Psi_y(v_j)\gets \psi(v_j+\alpha\iu), \qquad
\Psi_z(v_j)\gets \sigma(\iu v_j-\alpha)\psi(v_j+\alpha\iu)
\]
\State Precompute recovery terms on the grid:
\[
\ddot Y(x_i)\gets e^{x_i}\psi(-\iu),\qquad
\ddot Z(x_i)\gets -\frac{\eta\Delta t\,e^{x_i}\psi(-\iu)+\iu e^{x_i}\psi'(-\iu)}{\sigma\Delta t}
\]

\For{$k \gets n-1, n-2, \dots, 0$}
    \State \textbf{(Shift)} Compute boundary slopes (one-sided differences):
    \[
      y'_0 \gets \frac{-3Y_{k+1,0}^\Delta + 4Y_{k+1,1}^\Delta - Y_{k+1,2}^\Delta}{2\Delta x},\quad
      y'_{N-1} \gets \frac{3Y_{k+1,N-1}^\Delta - 4Y_{k+1,N-2}^\Delta + Y_{k+1,N-3}^\Delta}{2\Delta x}
    \]
    \State Solve for shifting parameters $A,B$ (enforcing periodicity of $\tilde u$ and $\tilde u'$ on $[x_0,x_{N-1}]$)
    \State Form the damped--shifted vector on the grid:
    \[
      \tilde Y_{k+1,i}^{(\alpha)} \gets e^{\alpha x_i}\bigl(Y_{k+1,i}^\Delta - (A e^{x_i}+B)\bigr), \quad i=0,\dots,N-1
    \]

    \State \textbf{(FFT step)} Compute centered DFT/IDFT using the phase factors $(-1)^i$:
    \Statex \hspace{\algorithmicindent} $ \widehat{\tilde Y}\gets \mathfrak{D}\bigl(((-1)^i\tilde Y_{k+1,i}^{(\alpha)})_{i=0}^{N-1}\bigr)$
    \Statex \hspace{\algorithmicindent} $ \widehat{\dot Y}\gets \widehat{\tilde Y}\odot \Psi_y$ \Comment{pointwise product}
    \Statex \hspace{\algorithmicindent} $ \widehat{\dot Z}\gets \widehat{\tilde Y}\odot \Psi_z$
    \Statex \hspace{\algorithmicindent} $ \dot Y_{k}^\Delta \gets ((-1)^i\,\mathfrak{D}^{-1}(\widehat{\dot Y}))_{i=0}^{N-1}$
    \Statex \hspace{\algorithmicindent} $ \dot Z_{k}^\Delta \gets ((-1)^i\,\mathfrak{D}^{-1}(\widehat{\dot Z}))_{i=0}^{N-1}$

    \State \textbf{(Undamp/unshift)} Recover $(\hat Y_k^\Delta,\bar Z_k^\Delta)$ on the grid:
    \[
      \hat Y_{k,i}^\Delta \gets e^{-\alpha x_i}\dot Y_{k,i}^\Delta + A\,\ddot Y(x_i) + B,
      \qquad
      \bar Z_{k,i}^\Delta \gets e^{-\alpha x_i}\dot Z_{k,i}^\Delta + A\,\ddot Z(x_i)
    \]

    \State \textbf{(Driver update)} Explicit Euler step for $Y$:
    \[
      \bar Y_{k,i}^\Delta \gets \hat Y_{k,i}^\Delta + f(X_{k,i}^\Delta,\hat Y_{k,i}^\Delta,\bar Z_{k,i}^\Delta)\Delta t
    \]
    \State Optional constraint (e.g.\ call payoff): $\bar Y_{k,i}^\Delta \gets \max(\bar Y_{k,i}^\Delta,0)$
    \State Set $Y_{k,i}^\Delta \gets \bar Y_{k,i}^\Delta$, $Z_{k,i}^\Delta \gets \bar Z_{k,i}^\Delta$ for all $i$
\EndFor

\State \Return $\{\bar Y_k^\Delta(x_i),\bar Z_k^\Delta(x_i)\}$
\end{algorithmic}
\end{algorithm}

	In Algorithm \ref{Algo2}, we fix the damping parameter and only update the shifting parameter at every time step which gives fast and efficient calculation. Comparing the results to those obtained using the convolution algorithm given by \citet{hyndman2017convolution}, where both the shifting and the damping parameters are adaptively updated at every time step, we find that the changes embodied in Algorithm~\ref{Algo2} increase stability.

	\subsection{Error analysis}
	We denote the convolution estimation for $u(x_n) = \mathbb{E}\left[Y\left|X_{t_{n}}=x_n\right.\right]$ 
	\begin{equation*}
		\bar u(x_n) \coloneqq (-1)^{n} \mathfrak{D}^{-1}\left[\left\{\mathfrak{D}\left[\left\{w_n(-1)^{n}Y^{\Delta}_i\right\}_{i=0}^{N-1}\right](v_j)\Psi\left(v_j\right)\right\}_{j=0}^{N-1}\right]_{n}.
	\end{equation*}
	Following Proposition 2.1 of \citet{wang_hynd:heston} and the fact that the characteristic function of Gaussian density decays as $\exp\left(-\frac{1}{2}\Sigma x^2 \right)$, we obtain the following lemma.
	\begin{lemma}[Error of the convolution method]\label{lem_conv}
		Assuming the integrable function $f(x)$ is bounded by $\bar f$ on $[-\frac{L}{2},\frac{L}{2}]$ and admits the Fourier expansion
		\begin{equation*}
			f(x) = \sum_{j=-\infty}^{\infty} F_j e^{-\iu j \frac{2\pi x}{L}},
		\end{equation*}
		with coefficients defined by
		\begin{equation*}
			F_j = \frac{1}{L}\int_{-\frac{L}{2}}^{\frac{L}{2}} f(x) e^{\iu j \frac{2\pi x}{L}} dx.
		\end{equation*}
		Suppose the discrete Fourier coefficient
		\begin{equation*}
			\bar{F}_j = \frac{\Delta x}{L}\sum_{k=0}^{N-1} f(x_k) e^{\iu j \frac{2\pi x_k}{L}},
		\end{equation*}
		has bounded error $\left|F_j - \bar F_j\right|\leq \epsilon_L N^{-m}$ for $m\geq 2$ and some constant $\epsilon_L>0$. Then the convolution method has an estimation error bounded by
		\begin{equation}\label{plain_err}
			\left|e_y\right|\coloneqq \left|u - \bar u\right| \leq C_L\left(\bar f\left(1 - \mathrm{erf}\left(\frac{\sigma\left(N-2\right)\pi}{L}\sqrt{\frac{\Delta t}{2}}\right)\right) + \epsilon L N^{-m}\right),
		\end{equation}
		for some constant $\epsilon>0$ depending on $L$ and $C_L=\frac{L}{\sigma\sqrt{2\pi\Delta t}}$ on the truncation regions $[-\frac{L}{2}, \frac{L}{2}]$ with discretization $\Delta = (\Delta t,\Delta x, \Delta v)$. The error function $\mathrm{erf}(x)$ is defined as
		\begin{equation*}
			\mathrm{erf}(x)\coloneqq \frac{2}{\sqrt{\pi}} \int_0^x e^{-u^2}du.
		\end{equation*}
	\end{lemma}
	\begin{remark}[Boundary problem]
		In Lemma \ref{lem_conv}, we use the bound  $\bar f$ to estimate the Fourier coefficients as the target function may be unbounded with increasing derivatives. The terminal function used in option pricing is a such function that is not in the Schwartz space of all functions whose derivatives are rapidly decreasing. The rapidly increasing derivatives of the log-underlying makes its Fourier coefficient $\left|F_j\right|$ increase as $x\rightarrow\infty$. To analyze the Fourier transform of such target function, the bound $\bar f$ could heavily influence the local error in addition to the truncation factor $L$ and the discretization factor $\Delta$. Similar to \citet{hyndman2017convolution} we will prove local and global error estimates, stability and convergence criteria, and illustrate boundary problem numerically in the application section.
	\end{remark}
	Lemma \ref{lem_conv} provides a general error estimate for the convolution method, and as we can see from (\ref{plain_err}) that if the target function is unbounded at one side, the value of $\bar f$ can be large and the convolution method yields poor approximation results at boundaries. This result can be viewed from the numerical results provided by \citet{hyndman2017convolution}. In the convolution method proposed by \citet{hyndman2017convolution}, a linear function was applied to modify the target function, which provides reliable accuracy in the center area of the truncation region but does not improve  performance on boundaries. That is, the error on the boundaries are becoming unbounded since the target function itself is an exponential-type function and the linear function is not compatible to shift  it as a bounded function. Considering that the shifted function yields  smoothly connected boundaries, we could find a smaller $\hat f< \bar f$ to bound the shifted function $\tilde f$. Following Lemma \ref{lem_conv} and consider the Fourier transform with damping and shifting schemes, the error estimate is then given by
	\begin{equation}
		\left|\tilde e_y\right|\leq C_{\alpha,L}e^{-\alpha x}\left(\hat f\left(1 - erf\left(\frac{\sigma\left(N-2\right)\pi}{L}\sqrt{\frac{\Delta t}{2}}\right)\right) + \epsilon L N^{-m}\right),\label{controlled_y_error}
	\end{equation}
	for the constant $C_{\alpha, L}=\frac{Le^{-\alpha\Delta t\left(\eta -\frac{1}{2}\sigma^2\alpha\right)}}{\sigma\sqrt{2\pi\Delta t}}$ and some value $\hat f \geq \left|\tilde f(x)\right|$ for all $x$. Then, the error for $Z$ is bounded by
	\begin{equation}
		\left|\tilde e_z\right| \leq C_{\alpha, L}\frac{\sigma e^{-\alpha x}}{|\alpha|}\left(\hat f\left(1 - erf\left(\frac{\sigma\left(N-2\right)\pi}{L}\sqrt{\frac{\Delta t}{2}}\right)\right) + \frac{1}{\sigma|\alpha| \sqrt{2\pi\Delta t}}\right) + \mathcal{O}\left(e^{-K\frac{N^2}{L^2}\Delta t}\right),\label{controlled_z_error}
	\end{equation}
	for the constant $K=\frac{\sigma^2\pi^2}{2}$.
	\begin{remark}[Error transfer with damping parameter]
		From (\ref{controlled_y_error}) and (\ref{controlled_z_error}) we observe that the shifting term reduces the truncation error given that $\hat f < \bar f$ while the damping term $e^{-\alpha x}$ has a side effect which makes the error increase as $x$ approaches to the right boundary for $\alpha<0$. The proper shifting function $h(x)$ is chosen to be similar to the terminal condition $g(x)$ of $Y_t(x)$. The shifting result gives us a periodic function which yields  smaller error estimates than using $\hat f$ to bound the Fourier coefficients, see Theorem 4.4 of \citet{vretblad2003fourier}, where the convolution error is of order $\mathcal{O}\left(N^{-1}\right)$.  We can rewrite (\ref{controlled_y_error}) and (\ref{controlled_z_error}) as
		\begin{align}
			\left|\tilde e_y\right| \leq& \mathcal{O}\left(N^{-1}\right),\label{controlled_y_error1}\\
			\left|\tilde e_z\right| \leq& \mathcal{O}\left(\Delta t^{-1}\right) + \mathcal{O}\left(N^{-1}\right) + \mathcal{O}\left(e^{-K\frac{N^2}{L^2}\Delta t}\right).\label{controlled_z_error1}
		\end{align}
	\end{remark}
	By Lemma \ref{lem_conv} and the Lipschitz condition of $f$, we obtain the local error estimation for the convolution method with the damping and shifting scheme.
	\begin{lemma}[Local error of the convolution method with damping and shifting]
		Suppose Assumption \ref{assumption1} is satisfied, the damping and the shifting schemes admits the following error estimation on the discretized region $\Delta = \left(\Delta t, \Delta x, \Delta v\right)$
		\begin{equation*}
			\left|Y(x) - \bar Y^\Delta(x)\right|\leq C_\alpha e^{-\alpha x} \left(\Delta x + \sqrt{\Delta t}\right),
		\end{equation*}
		for some constant $C_\alpha>0$ depending only on $\alpha$.
	\end{lemma}
	Let $\frac{\sigma\left(N-2\right)\pi}{L}\sqrt{\frac{\Delta t}{2}}$ be large enough and $erf\left(\frac{\sigma\left(N-2\right)\pi}{L}\sqrt{\frac{\Delta t}{2}}\right) \rightarrow 1$ such that the truncation error converges in (\ref{plain_err}), which yields the following convergence condition.
	\begin{proposition}[Stability and convergence]
		If the discretization $N$ and the truncation $L$ satisfies
		\begin{equation}\label{convergence_condition}
			N\geq \frac{L}{\sigma}\sqrt{\frac{2}{\Delta t}},
		\end{equation}
		then the convolution method is stable and convergent.
	\end{proposition}
	
	Next, we investigate the global error estimation and summarize it by the following Theorem.
	\begin{theorem}[Global error bounds]\label{thm_conv_bsde}
		Suppose $f\in\mathbb{C}^{1,2,2,2}\rightarrow\mathbb{R}$ satisfies Assumption~\ref{assumption1} and  (\ref{convergence_condition}) holds Then, the error of the convolution-FFT numerical solution of (\ref{FSDE})-(\ref{BSDE}) on the discretized region $\Delta = \left(\Delta t, \Delta x, \Delta v\right)$ is of the form
		\begin{equation*}
			\mathrm{err}(\Delta)\coloneqq \max_{0\leq i\leq n} \mathbb{E}\left[\sup_{t_i\leq t\leq t_{i+1}}\left|Y_t - \bar  Y^\Delta_{t_i}\right|\right]
			\leq \mathcal{O}\left(\Delta t^{\frac{1}{2}}\right) + \mathcal{O}\left(\Delta t\right) +  \mathcal{O}\left(\Delta t \Delta x\right) + \mathcal{O}\left(\Delta t e^{-K\frac{\Delta t}{\Delta x^2}}\right),
		\end{equation*}
		for the constant $K=\frac{\sigma^2\pi^2}{2}$.
	\end{theorem}
	\noindent Proof: see Appendix \ref{pf_thm_conv_bsde}.

Theorem~\ref{thm_conv_bsde} decomposes the global error into three contributions: time-discretization error of order $\mathcal{O}(\Delta t^{1/2})+\mathcal{O}(\Delta t)$, spatial discretization error of order $\mathcal{O}(\Delta t\,\Delta x)$, and a truncation term of order $\mathcal{O}\!\left(\Delta t\,\exp{(-K\Delta t/\Delta x^{2})}\right)$. The boundary-control strategy developed in Section~\ref{sec3.1} is designed to prevent the truncation and periodic-extension effects from dominating the numerical solution when the target functions are non-periodic and unbounded on the full line.

Compared with the adaptive damping and shifting strategy in \citet{hyndman2017convolution}, we fix the damping parameter and update only the shifting parameters at each time step. This avoids step-to-step variation in the damping exponent, which can amplify boundary artifacts in practice. Moreover, the truncation term in Theorem~\ref{thm_conv_bsde} decays rapidly as the spatial mesh is refined under the stability condition \eqref{convergence_condition}, and the numerical experiments in the next section confirm that the resulting scheme achieves improved accuracy on comparable grids.

	\section{Numerical result of option pricing}\label{sec:application}
\begin{figure}[H]
		\centering
		\caption{Call option price and delta errors - convolution method of \citet{hyndman2017convolution}}
		\vspace*{-0.4cm}
		\includegraphics[width=1\linewidth]{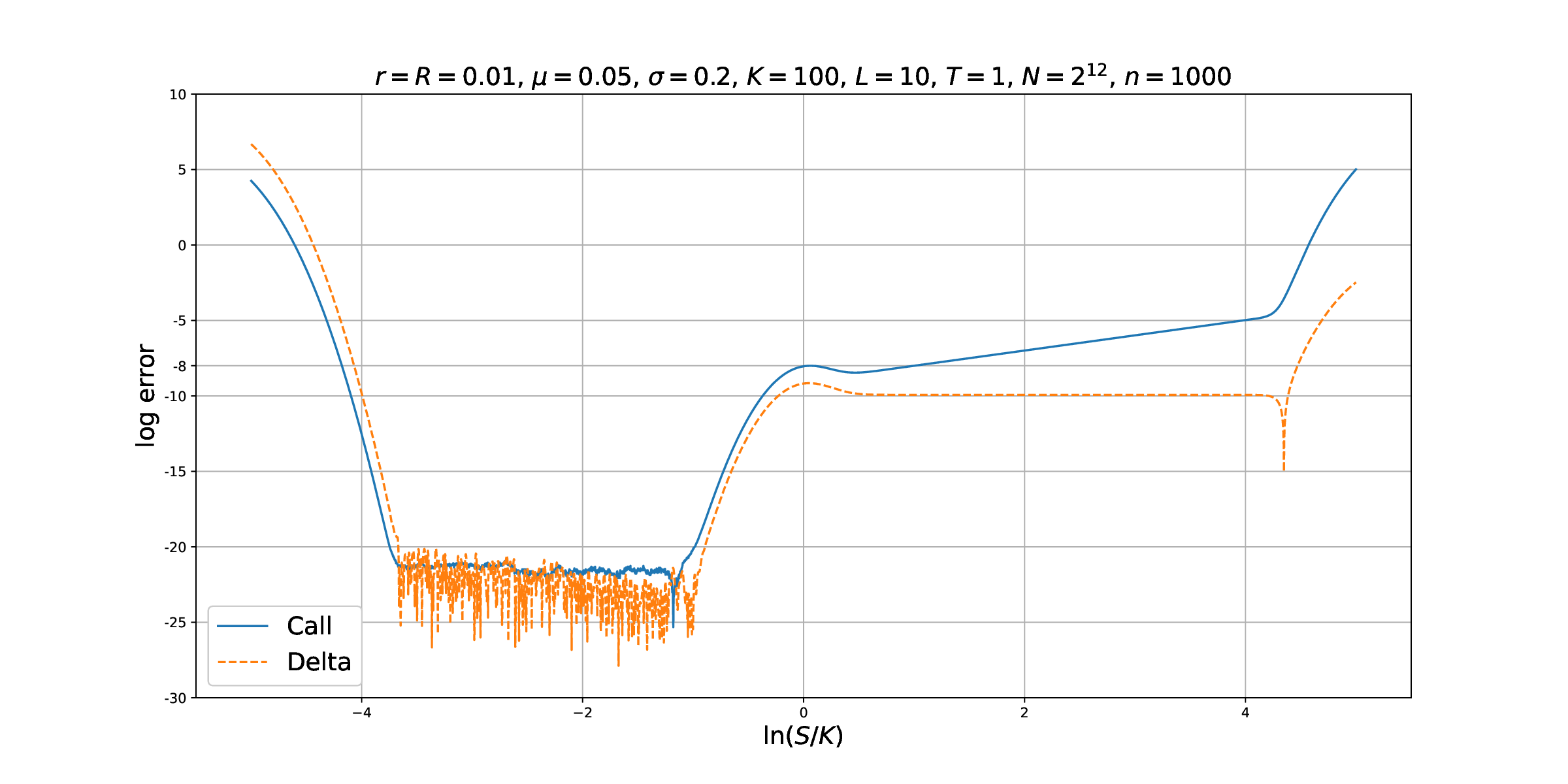}
		\vspace*{-1cm} 
		\label{pic1} %
	\end{figure}

	Suppose the underlying asset $S_t$ pays constant dividend with constant $\mu$, $d$, $\sigma$ which is defined as follows
	\begin{equation}
		{S}_t = {S}_0 \exp\left\{\left(\mu - d - \frac{1}{2}\text{tr}\left(\sigma\sigma'\right)\right)t + \sigma W_t\right\}.
	\end{equation}
	The corresponding logarithm of the stock prices $X_t = \ln S_t$ is defined in (\ref{FSDE}) with $\eta_t = \mu - d - \frac{1}{2}\text{tr}\left(\sigma\sigma'\right)$ and $X_0=\ln S_0$. In the market with borrowing rate $R_t$, the driver for the American call option price is given by
	\begin{equation}
		f(t,{x}, y,{z}) = -r_t y - \left(\mu_t - d_t -  r\right){\sigma_t}^{-1} {z} + \left(R_t - r_t\right)\left(y- \text{tr}\left(\sigma_t^{-1}{z}\right)\right)^-.
	\end{equation}
	The terminal condition of options in European or American type is
	\begin{equation*}
		g(x) = \left(e^{x} - K\right)^+.
	\end{equation*}
	We choose that $S_0=100$, $K=100$, $d=0$, $r=R=0.01$, $\mu=0.05$, $\sigma=0.2$, and $T=1$. In Figure \ref{pic1}, we replicate the result and the method provided by \citet{hyndman2017convolution}. We construct the mesh of the spatial variable $X$ with $N$ discretized values on $X_0 + [-\frac{L}{2},\frac{L}{2}]$ and the backward iterations over the time with $n$ steps for $\Delta t=\frac{T}{n}$.  Since $d=0$ and $R=r$, our numerical results specialize to the Black-Scholes European call option price.  
	
	\begin{figure}[H]
		\centering
		\caption{Call option price and delta errors - convolution method with boundary error control.}
		\vspace*{-0.4cm} 
		\includegraphics[width=1\linewidth]{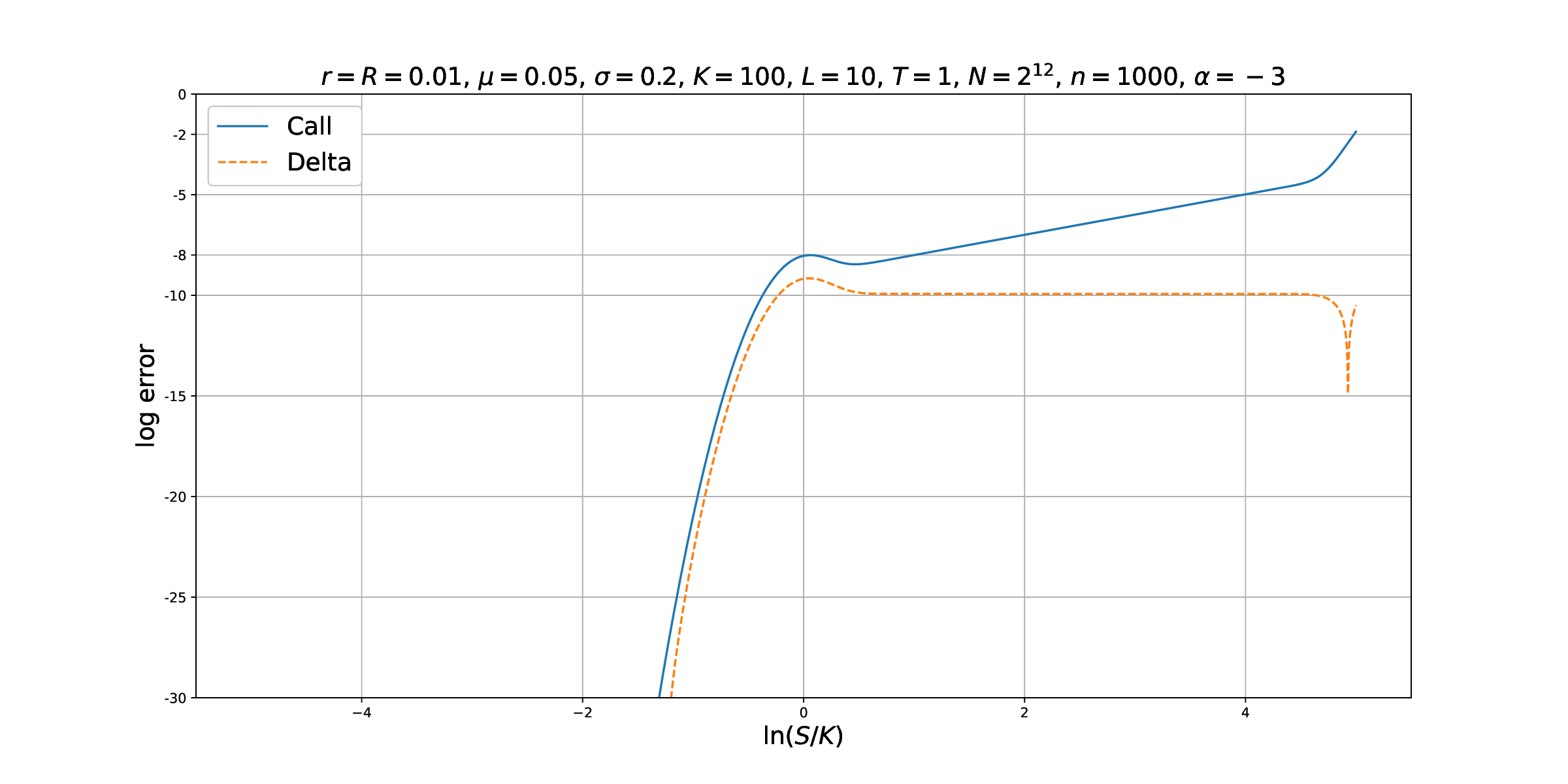}
		\vspace*{-1cm} 
		\label{pic2} %
	\end{figure}

	As we can see from Figure \ref{pic1} and Figure \ref{pic2}, the CFFT--BSDE method with boundary control substantially decreases errors at the boundaries, effectively eliminating boundary error for deeply out-of-the-money (OTM) options. The  CFFT--BSDE method with boundary control also provides more stable results compared to \citet{hyndman2017convolution} method which has wide range of damped oscillation shown in Figure \ref{pic1}.  Table~\ref{delta_sweep_ATM} examines the impacts of different parameters on the hedging ratio.  We calculate the hedge ratios using $\Delta = Z/(\sigma S)$ based on (\ref{markovian_repZ}) and using finite differences derived from the CFFT-BSDE approximation of $Y_{t}$.
	
	\begin{figure}[H]
		\centering
		\caption{Call option delta surface - convolution method with boundary error control}
		\includegraphics[width=0.5\linewidth]{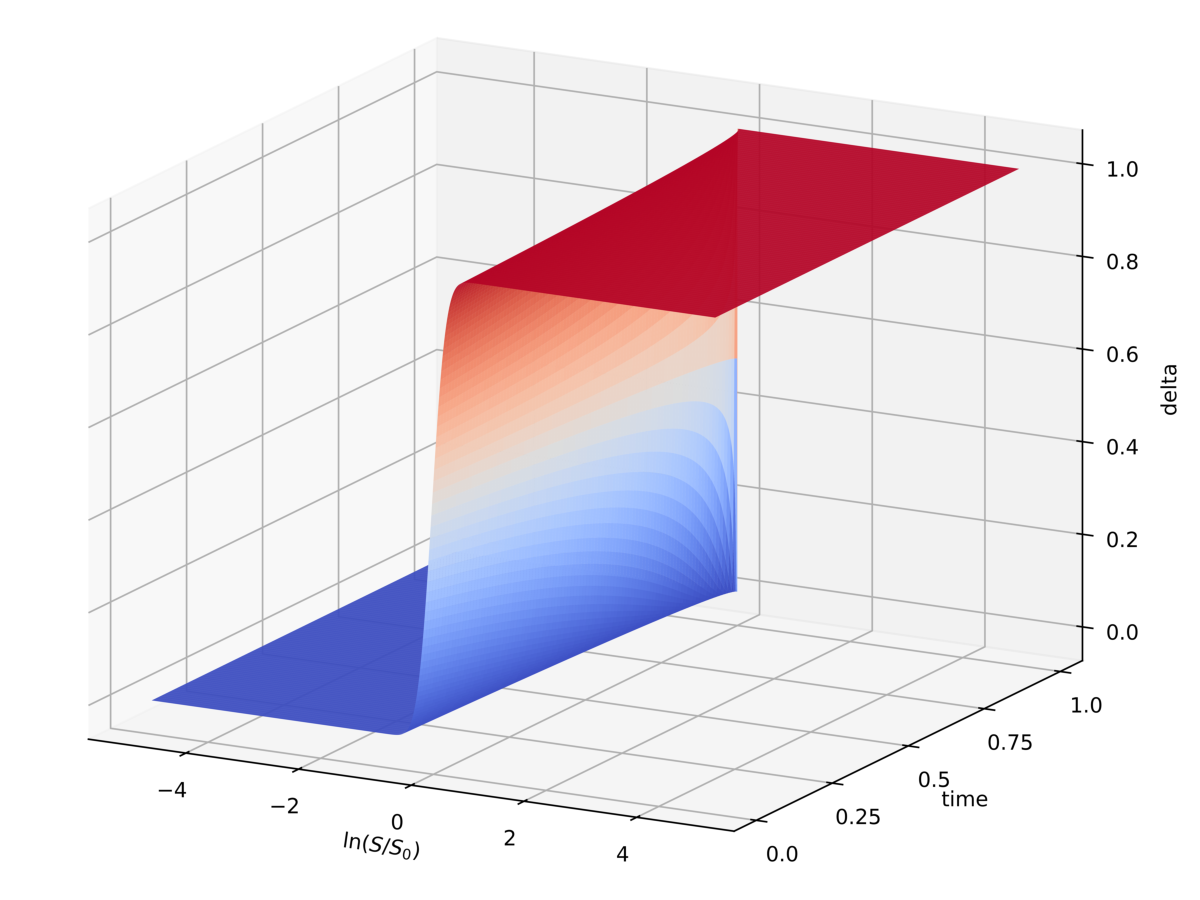}
		\label{pic3} %
	\end{figure}

	The delta surface provided in \citet{hyndman2017convolution} does not show the full picture on the truncation region since the value of delta explodes at boundary as well. Figure \ref{pic3} demonstrates that our method displays the entire delta surface which is smooth and very accurate at the boundaries.
	
	\begin{table}[!t]
\centering
\small
\setlength{\tabcolsep}{4.5pt}
\renewcommand{\arraystretch}{1.15}
  \caption{ATM Delta errors for the CFFT/BSDE method (via $Z$ and finite differences) against Black--Scholes.
  $S_{0}=100$, $K=100$, $r=0.01$, $\sigma=0.2$, $T=1$. }
\label{tab:delta_sweep}
\begin{tabular}{cccccccccc}
\toprule
$n$ & $L$ & $N$ & $\Delta_{\mathrm{BS}}$ & $\Delta_{Z}$ & $|\Delta_{Z}-\Delta_{\mathrm{BS}}|$ & $\frac{|\Delta_{Z}-\Delta_{\mathrm{BS}}|}{|\Delta_{\mathrm{BS}}|}$ & $\Delta_{\mathrm{FD}}$ & $|\Delta_{\mathrm{FD}}-\Delta_{\mathrm{BS}}|$ & $\frac{|\Delta_{\mathrm{FD}}-\Delta_{\mathrm{BS}}|}{|\Delta_{\mathrm{BS}}|}$ \\
\midrule
1000 & 10.000000 & $2^{10}$ & 0.559618 & 0.559619 & 1.153e-06 & 2.061e-06 & 0.559613 & 4.418e-06 & 7.894e-06 \\
1000 & 10.000000 & $2^{11}$ & 0.559618 & 0.559622 & 4.616e-06 & 8.249e-06 & 0.559617 & 9.785e-07 & 1.749e-06 \\
1000 & 10.000000 & $2^{12}$ & 0.559618 & 0.559623 & 5.351e-06 & 9.562e-06 & 0.559617 & 2.448e-07 & 4.375e-07 \\
1000 & 12.000000 & $2^{10}$ & 0.559618 & 0.559616 & 1.266e-06 & 2.263e-06 & 0.559611 & 6.810e-06 & 1.217e-05 \\
1000 & 12.000000 & $2^{11}$ & 0.559618 & 0.559622 & 4.178e-06 & 7.466e-06 & 0.559616 & 1.415e-06 & 2.528e-06 \\
1000 & 12.000000 & $2^{12}$ & 0.559618 & 0.559623 & 5.244e-06 & 9.370e-06 & 0.559617 & 3.525e-07 & 6.299e-07 \\
1000 & 14.000000 & $2^{10}$ & 0.559618 & 0.559613 & 4.369e-06 & 7.807e-06 & 0.559608 & 9.868e-06 & 1.763e-05 \\
1000 & 14.000000 & $2^{11}$ & 0.559618 & 0.559621 & 3.634e-06 & 6.494e-06 & 0.559616 & 1.956e-06 & 3.495e-06 \\
1000 & 14.000000 & $2^{12}$ & 0.559618 & 0.559623 & 5.116e-06 & 9.142e-06 & 0.559617 & 4.797e-07 & 8.571e-07 \\
2000 & 10.000000 & $2^{10}$ & 0.559618 & 0.559616 & 1.925e-06 & 3.441e-06 & 0.559613 & 4.698e-06 & 8.396e-06 \\
2000 & 10.000000 & $2^{11}$ & 0.559618 & 0.559619 & 1.801e-06 & 3.218e-06 & 0.559617 & 9.958e-07 & 1.779e-06 \\
2000 & 10.000000 & $2^{12}$ & 0.559618 & 0.559620 & 2.553e-06 & 4.562e-06 & 0.559617 & 2.448e-07 & 4.375e-07 \\
2000 & 12.000000 & $2^{10}$ & 0.559618 & 0.559613 & 4.467e-06 & 7.982e-06 & 0.559610 & 7.213e-06 & 1.289e-05 \\
2000 & 12.000000 & $2^{11}$ & 0.559618 & 0.559619 & 1.310e-06 & 2.341e-06 & 0.559616 & 1.485e-06 & 2.653e-06 \\
2000 & 12.000000 & $2^{12}$ & 0.559618 & 0.559620 & 2.445e-06 & 4.370e-06 & 0.559617 & 3.525e-07 & 6.299e-07 \\
2000 & 14.000000 & $2^{10}$ & 0.559618 & 0.559610 & 7.639e-06 & 1.365e-05 & 0.559607 & 1.034e-05 & 1.848e-05 \\
2000 & 14.000000 & $2^{11}$ & 0.559618 & 0.559618 & 6.851e-07 & 1.224e-06 & 0.559616 & 2.107e-06 & 3.765e-06 \\
2000 & 14.000000 & $2^{12}$ & 0.559618 & 0.559620 & 2.318e-06 & 4.142e-06 & 0.559617 & 4.798e-07 & 8.574e-07 \\
5000 & 10.000000 & $2^{10}$ & 0.559618 & 0.559614 & 3.791e-06 & 6.775e-06 & 0.559613 & 4.885e-06 & 8.730e-06 \\
5000 & 10.000000 & $2^{11}$ & 0.559618 & 0.559618 & 7.899e-08 & 1.412e-07 & 0.559617 & 1.039e-06 & 1.856e-06 \\
5000 & 10.000000 & $2^{12}$ & 0.559618 & 0.559619 & 8.740e-07 & 1.562e-06 & 0.559617 & 2.451e-07 & 4.380e-07 \\
5000 & 12.000000 & $2^{10}$ & 0.559618 & 0.559611 & 6.346e-06 & 1.134e-05 & 0.559610 & 7.413e-06 & 1.325e-05 \\
5000 & 12.000000 & $2^{11}$ & 0.559618 & 0.559617 & 4.101e-07 & 7.328e-07 & 0.559616 & 1.526e-06 & 2.727e-06 \\
5000 & 12.000000 & $2^{12}$ & 0.559618 & 0.559618 & 7.638e-07 & 1.365e-06 & 0.559617 & 3.553e-07 & 6.348e-07 \\
5000 & 14.000000 & $2^{10}$ & 0.559618 & 0.559608 & 9.524e-06 & 1.702e-05 & 0.559607 & 1.055e-05 & 1.885e-05 \\
5000 & 14.000000 & $2^{11}$ & 0.559618 & 0.559617 & 1.062e-06 & 1.898e-06 & 0.559616 & 2.175e-06 & 3.887e-06 \\
5000 & 14.000000 & $2^{12}$ & 0.559618 & 0.559618 & 6.281e-07 & 1.122e-06 & 0.559617 & 4.908e-07 & 8.770e-07 \\
\bottomrule
\end{tabular} \label{delta_sweep_ATM}
\end{table}

	\section{Conclusion}\label{sec:conclusion}
	
	In this paper, we propose a numerical method that improves the performance using convolution method in solving BSDEs and demonstrate that the boundary error is significantly reduced by our method. This numerical method provides a new approach to improve the accuracy of convolution method with the fast Fourier transform. Our motivation is focused on the transformation of the target function using a shifting function which is similar to the terminal of the BSDEs and is able to map the target function as a bounded periodic function. We have studied the application of the convolution method in valuing options through the framework of BSDEs and provided detailed error analysis including three parts from extrapolation, truncation to discretization. The numerical result shows that our method has better accuracy than the original method given by \citet{hyndman2017convolution}. Both the theoretical analysis and numerical result show us that the boundary error still increases with respect to the truncation domain, however, the boundary error is well controlled by using our method in the Fourier transform for the unbounded and non-periodic problem. Therefore, our method is feasible to apply on more general BSDEs problems. Future work will investigate extensions of this approach to more general BSDEs and to higher dimensional problems.

        \nocite{XGao-PHDthesis2021}
	\bibliographystyle{chicago}
	\bibliography{mybib,newbib}

\begin{thebibliography}{}

\bibitem[\protect\citeauthoryear{Antonelli}{Antonelli}{1993}]{antonelli1993backward}
Antonelli, F. (1993).
\newblock {\em Backward forward stochastic differential equations}.
\newblock Ph.\ D. thesis, Purdue University.

\bibitem[\protect\citeauthoryear{Bender and Denk}{Bender and
  Denk}{2007}]{bender2007forward}
Bender, C. and R.~Denk (2007).
\newblock A forward scheme for backward sdes.
\newblock {\em Stochastic processes and their applications\/}~{\em 117\/}(12),
  1793--1812.

\bibitem[\protect\citeauthoryear{Bouchard, Elie, and Touzi}{Bouchard
  et~al.}{2009}]{BouchardElieTouzi+2009+91+124}
Bouchard, B., R.~Elie, and N.~Touzi (2009).
\newblock Discrete-time approximation of bsdes and probabilistic schemes for
  fully nonlinear pdes.
\newblock In H.~Albrecher, W.~J. Runggaldier, and W.~Schachermayer (Eds.), {\em
  Advanced Financial Modelling}, pp.\  91--124. De Gruyter.

\bibitem[\protect\citeauthoryear{Bouchard and Touzi}{Bouchard and
  Touzi}{2004}]{bouchard2004discrete}
Bouchard, B. and N.~Touzi (2004).
\newblock Discrete-time approximation and monte-carlo simulation of backward
  stochastic differential equations.
\newblock {\em Stochastic Processes and Their Applications\/}~{\em 111\/}(2),
  175--206.

\bibitem[\protect\citeauthoryear{Carr and Madan}{Carr and
  Madan}{1999}]{carr1999option}
Carr, P. and D.~Madan (1999).
\newblock Option valuation using the fast {F}ourier transform.
\newblock {\em Journal of Computational Finance\/}~{\em 2\/}(4), 61--73.

\bibitem[\protect\citeauthoryear{Douglas, Ma, and Protter}{Douglas
  et~al.}{1996}]{douglas1996numerical}
Douglas, J., J.~Ma, and P.~Protter (1996).
\newblock Numerical methods for forward-backward stochastic differential
  equations.
\newblock {\em Annals of Applied Probability\/}~{\em 6\/}(3), 940--968.

\bibitem[\protect\citeauthoryear{Gao}{Gao}{2021}]{XGao-PHDthesis2021}
Gao, X. (2021, March).
\newblock {\em Stochastic control, numerical methods, and machine learning in
  finance and insurance}.
\newblock Ph.\ D. thesis, Concordia University.

\bibitem[\protect\citeauthoryear{Gao and Hyndman}{Gao and
  Hyndman}{2025}]{wang_hynd:heston}
Gao, X. and C.~B. Hyndman (2025).
\newblock Convolution-{FFT} for option pricing in the {H}eston model.
\newblock arXiv preprint.
\newblock 2512.05326.

\bibitem[\protect\citeauthoryear{Huijskens, Ruijter, and Oosterlee}{Huijskens
  et~al.}{2016}]{huijskens2016efficient}
Huijskens, T., M.~J. Ruijter, and C.~W. Oosterlee (2016).
\newblock Efficient numerical {F}ourier methods for coupled forward--backward
  {SDEs}.
\newblock {\em Journal of Computational and Applied Mathematics\/}~{\em 296},
  593--612.

\bibitem[\protect\citeauthoryear{Hyndman and Oyono~Ngou}{Hyndman and
  Oyono~Ngou}{2017}]{hyndman2017convolution}
Hyndman, C.~B. and P.~Oyono~Ngou (2017).
\newblock A convolution method for numerical solution of backward stochastic
  differential equations.
\newblock {\em Methodology and Computing in Applied Probability\/}~{\em
  19\/}(1), 1--29.

\bibitem[\protect\citeauthoryear{Kobylanski}{Kobylanski}{2000}]{kobylanski2000backward}
Kobylanski, M. (2000).
\newblock Backward stochastic differential equations and partial differential
  equations with quadratic growth.
\newblock {\em Annals of Probability\/}~{\em 28\/}(2), 558--602.

\bibitem[\protect\citeauthoryear{Lemor, Gobet, and Warin}{Lemor
  et~al.}{2006}]{lemor2006rate}
Lemor, J.-P., E.~Gobet, and X.~Warin (2006).
\newblock Rate of convergence of an empirical regression method for solving
  generalized backward stochastic differential equations.
\newblock {\em Bernoulli\/}~{\em 12\/}(5), 889--916.

\bibitem[\protect\citeauthoryear{Lord and Kahl}{Lord and
  Kahl}{2006}]{lord2007optimal}
Lord, R. and C.~Kahl (2006, 01).
\newblock Optimal {F}ourier inversion in semi-analytical option pricing.
\newblock {\em Tinbergen Institute, Tinbergen Institute Discussion
  Papers\/}~{\em 10}, 1--23.

\bibitem[\protect\citeauthoryear{Ma, Morel, and Yong}{Ma
  et~al.}{1999}]{ma1999forward}
Ma, J., J.-M. Morel, and J.~Yong (1999).
\newblock {\em Forward-backward stochastic differential equations and their
  applications}.
\newblock Number 1702 in Lecture Notes in Mathematics. Berlin Heidelberg:
  Springer Verlag.

\bibitem[\protect\citeauthoryear{Oyono~Ngou and Hyndman}{Oyono~Ngou and
  Hyndman}{2022}]{polynice2022convolution}
Oyono~Ngou, P. and C.~B. Hyndman (2022).
\newblock A {F}ourier interpolation method for numerical solution of {FBSDE}s:
  {G}lobal convergence, stability, and higher order discretizations.
\newblock {\em Journal of Risk and Financial Management\/}~{\em 15\/}(9), 388.

\bibitem[\protect\citeauthoryear{Pardoux and Peng}{Pardoux and
  Peng}{1990}]{pardoux1990adapted}
Pardoux, E. and S.~Peng (1990).
\newblock Adapted solution of a backward stochastic differential equation.
\newblock {\em Systems \& Control Letters\/}~{\em 14\/}(1), 55--61.

\bibitem[\protect\citeauthoryear{Risken}{Risken}{1996}]{risken1996fokker}
Risken, H. (1996).
\newblock {\em The {F}okker-{P}lanck Equation\/} (2nd ed.).
\newblock Berlin Heidelberg: Springer.

\bibitem[\protect\citeauthoryear{Vretblad}{Vretblad}{2003}]{vretblad2003fourier}
Vretblad, A. (2003).
\newblock {\em {F}ourier {A}nalysis and its {A}pplications}, Volume 223.
\newblock New York: Springer.

\bibitem[\protect\citeauthoryear{Yong}{Yong}{1997}]{yong1997finding}
Yong, J. (1997).
\newblock Finding adapted solutions of forward-backward stochastic differential
  equations: method of continuation.
\newblock {\em Probability Theory and Related Fields\/}~{\em 107\/}(4),
  537--572.

\bibitem[\protect\citeauthoryear{Zhang}{Zhang}{2004}]{zhang2004numerical}
Zhang, J. (2004).
\newblock A numerical scheme for {BSDEs}.
\newblock {\em Annals of Applied Probability\/}~{\em 14\/}(1), 459--488.

\end{thebibliography}
	
	\appendix
	
	\section{Appendix}
	We provide the technical results for the proofs.
	\subsection{Proof of Theorem \ref{thm_conv_bsde}}
	\begin{proof}\label{pf_thm_conv_bsde}
		According to \citet{BouchardElieTouzi+2009+91+124}, the explicit Euler scheme 
		\begin{equation*}
			\left\{
			\begin{aligned}
				Z_{t_i} =& \frac{1}{\Delta t}\mathbb{E}\left[Y_{i+1} \Delta {W}_{t_i}\right],\\
				Y_{t_i} =& \mathbb{E}\left[Y_{t_{i+1}}\right] + \Delta t f\left({t_i},X_{t_i},\mathbb{E}\left[Y_{t_{i+1}}\right],{Z}_{t_i}\right),
			\end{aligned}
			\right.
		\end{equation*}
		admits an error in order $\mathcal{O}\left(\Delta t^{\frac{1}{2}}\right)$ for any $t\in[t_i,t_{i+1}]$
		\begin{equation}\label{error_Euler}
			\mathbb{E}\left[\sup_{t_i\leq t\leq t_{i+1}}\Big|Y_t - Y_{t_i}\Big|\right] \leq C\left(\Delta t^{\frac{1}{2}}\right).
		\end{equation}
		Since $f$ is Lipschitz continuous, we obtain that
		\begin{align}\label{error_conv}
			\left|Y_{t_i} - Y^\Delta_{t_i}\right| \leq& C\Delta t\left(\left|\mathbb{E}\left[Y_{t_{i+1}}\right] - \dot{\mathbb{E}}\left[Y_{t_{i+1}}\right]\right| + \left|Z_{t_{i+1}} - Z^\Delta_{t_{i+1}}\right|\right).
		\end{align}
		According to (\ref{controlled_y_error1}) and (\ref{controlled_z_error1}), we can rewrite (\ref{error_conv}) as
		\begin{align}\label{error_conv2}
			\left|Y_{t_i} - Y^\Delta_{t_i}\right| \leq& C\Delta t\left(N^{-1} + |\alpha| + \frac{1}{\sigma\sqrt{2\pi\Delta t}} + e^{-K\frac{N^2}{L^2}\Delta t}\right)\\
			\leq& C\left(\Delta t^{\frac{1}{2}} + \Delta t + \Delta t N^{-1} + \Delta t e^{-K\frac{N^2}{L^2}\Delta t}\right).
		\end{align}
		Combining (\ref{error_Euler}) and (\ref{error_conv2}), we obtain
		\begin{align}
			err(\Delta)\coloneqq& \max_{0\leq i\leq n} \mathbb{E}\left[\sup_{t_i\leq t\leq t_{i+1}}\left|Y_t - Y^\Delta_{t_i}\right|\right]\nonumber\\
			\leq& \max_{0\leq i\leq n} \mathbb{E}\left[\sup_{t_i\leq t\leq t_{i+1}}\Big|Y_t - Y_{t_i}\Big| + \left|Y_{t_i} - Y^\Delta_{t_i}\right|\right]\nonumber\\
			\leq& \mathcal{O}\left(\Delta t^{\frac{1}{2}}\right) + \mathcal{O}\left(\Delta t\right) +  \mathcal{O}\left(\Delta t N^{-1}\right) + \mathcal{O}\left(\Delta t e^{-K\frac{N^2}{L^2}\Delta t}\right)\nonumber\\
			\leq& \mathcal{O}\left(\Delta t^{\frac{1}{2}}\right) + \mathcal{O}\left(\Delta t\right) +  \mathcal{O}\left(\Delta t \Delta x\right) + \mathcal{O}\left(\Delta t e^{-K\frac{\Delta t}{\Delta x^2}}\right).\nonumber
		\end{align}

	\end{proof}
	
\end{document}